\documentclass[12pt,14paper]{article}
\usepackage{graphicx}
\usepackage{t1enc}
\usepackage[cmex10]{amsmath}
\usepackage{amssymb}
\usepackage{amsmath}
\usepackage{epsfig} 
\usepackage{cite}

\newtheorem{theorem}{Theorem}
\newtheorem{definition}{Definition}
\newtheorem{lemma}{Lemma}

\def\square{\hbox{\vrule\vbox{\hrule\phantom{o}\hrule}\vrule}}

\def\Re{\mathop{\rm Re}\nolimits}

\def\argmin{\mathop{\rm argmin}\nolimits}
\newcommand{\trn}{^{\scriptscriptstyle {\rm T}}}

\begin{document}
\title{
Decentralized adaptive synchronization in nonlinear dynamical networks with nonidentical nodes}
\author{Alexander~L.~Fradkov and~Ibragim~A.~Junussov
\thanks{A. L. Fradkov is with the Institute for Problems of Mechanical Engineering,
Russian Academy of Sciences,
61 Bolshoy ave V.O., St. Petersburg, 199178, Russia, fax: +7(812)321-4771, e-mail: fradkov@mail.ru}%
\thanks{I. A. Junussov is with Department of Theoretical Cybernetics, 
St.Petersburg State University, Universitetsky prospekt, 28, Peterhof, 198504, St. Petersburg, Russia,
e-mail: dxdtfxut@gmail.com}%
\thanks{The work was supported by the Russian Foundation for Basic Research (RFBR),
grants  07-01-92166, 08-01-00775; by Russian Academy of Sciences (RAS) youth
research program; by the program of basic research of OEMPPU RAS \#2 "Control and 
safety in energy and technical systems"; by the Council for grants of the RF President 
to support young Russian researchers and leading scientific schools (project NSh-2387.2008.1).}
}

\maketitle
\thispagestyle{empty}
\pagestyle{empty}

\setcounter{page}{1}
\begin{abstract}
\noindent For a network of interconnected nonlinear dynamical
systems an adaptive leader-follower output feedback synchronization 
problem is
considered. The proposed structure of decentralized controller and
adaptation algorithm is based on speed-gradient and
passivity. Sufficient conditions of
synchronization for nonidentical nodes are established. An example of synchronization of
the network of nonidentical Chua systems is analyzed. The main contribution of
the paper is adaptive controller design and analysis under conditions of
incomplete measurements, incomplete control and uncertainty.

\end{abstract}
\section{Introduction}
Adaptive synchronization of networked dynamical systems has 
attracted a growing interest during recent years 
\cite{LuChen05,YaoHillGuanWang06,ZhouLuLu06,ZhongDimirovskiZhao07}. It is 
motivated by a broad area of potential applications: formation control, 
cooperative control, control of power networks, communication networks,
production networks, etc. Existing works 
\cite{LuChen05,YaoHillGuanWang06,ZhouLuLu06,ZhongDimirovskiZhao07} 
and others are dealing with full state feedback and linear interconnections. 
The solutions are based on Lyapunov functions formed as sum of Lyapunov 
functions for local subsystems. As for adaptive control 
algorithms they are based on either local (decentralized
~\cite{Ioannou,Mirkin,GavelSiljak,WenSoh,JainKhorrami,Jiang,Fradkov_1990,FMN99}) or nearest
neighbor (described by an information graph~\cite{FaxMurray04,YoshiokaNamerikawa_Obs_based08,
ChopraSpong06,ChopraSpong08}) strategies.

Despite a great interest in control of network, only a restricted
class of them is currently solved. E.g. in existing papers
mainly linear  models of subsystems are considered
\cite{FaxMurray04, YoshiokaNamerikawa_Obs_based08}. 
In nonlinear case only passive or passifiable systems are studied
and control is organized according to information graph, i.e.
not completely decentralized~\cite{ChopraSpong06,ChopraSpong08}. 
Availability of the whole state
vector for measurement as well as appearance of control in all
equations for all nodes is assumed in decentralized stability and synchronization problems \cite
{LuChen05,YaoHillGuanWang06,ZhouLuLu06,ZhongDimirovskiZhao07}.
Powerful passivity based approaches are  not developed for adaptive synchronization problems.

In this paper we consider the problem of master-slave (leader-follower) synchronization in
a network of nonidentical systems in Lurie form where system models 
can be split into linear and nonlinear
parts. Case of identical nodes is studied in~\cite{FradkovJunussov_AIT08}.
Linearity of interconnections is not assumed; links between
subsystems can also be nonlinear. In the contrary to known works
on adaptive synchronization of networks, see
~\cite{ZhouLuLu06,ZhongDimirovskiZhao07}, only some output
function is available and control appears only in a part of the
system equations. It is also assumed that some plant parameters
are unknown.The leader subsystem is assumed to be isolated and
the control objective is to approach the trajectory of the
leader subsystem by all other ones under conditions of
uncertainty. Interconnection
functions are assumed to be Lipschitz continuous.

The results of~\cite{Fradkov_1990,FMN99} are employed to solve the posed problem. Adaptation
algorithm is designed by the speed-gradient method. It is shown
that the control goal is achieved under leader passivity
condition, if the interconnection strengths satisfy some inequalities.

The results are illustrated by example of synchronization 
in network of nonidentical Chua curcuits.

\section{Auxiliary results}
\subsection{Yakubovich-Kalman Lemma}
We need Yakubovich-Kalman Lemma in following form, see~\cite{YakubovichLeonovGelig}.
\begin{lemma}
\label{Yakub_lemma}
Let $A,B,C$ be $n\times n, n\times m, n\times m$ real matrices and $u\in\mathbb{R}^m,$ 
$\chi(s)=$\\
$C\trn (sI_n-A)^{-1}B, \mathrm{rank}\, B=m.$ 
Then the following statements are equivalent:

1) there exists matrix $H=H\trn>0$ such that
\begin{equation}
\label{HA_Yak}
HA+A\trn H<0, HB=C;
\end{equation}

2) polinomial $\det(sI_n-A)$ is Hurwitz and following frequency domain conditions hold
$$\mathrm{Re}\, u\trn \chi(i\omega)u>0,\quad \lim_{\omega\to\infty}\omega^2\,\mathrm{Re}\, u\trn\chi(i\omega)u>0$$
for all $\omega\in\mathbb{R}^1,$ $u\in\mathbb{R}^m, u\neq 0.$
\end{lemma}

\subsection{Speed gradient algorithm in decentralized control}
\label{section_th_2_18}
In order to present the main syncronization result of this paper we need to formulate problem 
statement of decentralized control and Theorem $2.18$ from~\cite{Fradkov_1990} which can also be derived
from Theorem $7.6$ from~\cite{FMN99}.

Consider\footnote{
    In this paper norms are Euclidean, $\mathrm{col}(x_1,\ldots,x_d)$ 
    stands for column vector with components consisting of components of $x_i, i=1,\ldots,d.$} 
a system $\mathcal{S}$ consisting of $d$ interconnected subsystems $\mathcal{S}_i,$ dynamics of 
each being described by the following equation:
\begin{equation}
\label{obj_equat_2_18}
\dot x_i = F_i(x_i, \tau_i, t) + h_i(x, \tau, t),\qquad i = 1,\ldots,d,
\end{equation} 
where $x_i\in \mathbb{R}^{n_i}$ -- state vector, $\tau_i\in\mathbb{R}^{m_i}$ - vector of inputs (tunable
parameters) of subsystem, $x = \mathrm{col}(x_1,\ldots,x_d)\in\mathbb{R}^n,$ 
$\tau = \mathrm{col}(\tau_1,\ldots,\tau_d)\in\mathbb{R}^m$
 - aggregate state and input vectors of system S, $n = \sum n_i,$ $m = \sum m_i.$ Vector-function $F_i(\cdot)$
describes local dynamics of subsystem $S_i,$  and vectors $h_i(\cdot)$ describe interconnection between subsystems.

Let $Q_i(x_i, t)\geq 0, i=1,\ldots,d$ be local goal functions and let the control goal be:
\begin{equation}
\label{goal_funcs_2_18}
\lim_{t\to\infty} Q_i(x_i, t) = 0,\qquad i = 1,\ldots, d.
\end{equation}
For all $i=1,\ldots,d$ we assume existence of smooth vector functions $x^*_i(t)$ such that $Q_i(x^*_i(t), t) \equiv 0,$ i.e.
$x^*_i = \argmin_{x_i} Q_i(x_i, t).$ Decentralized speed-gradient algorithm is introduced as follows:
\begin{equation}
\label{asg_2_18}
\dot\tau_i = -\Gamma_i\nabla_{\tau_i}\omega_i(x_i, \tau_i, t), \qquad i = 1,\ldots, d,
\end{equation}
where
$$
\omega_i(x_i, \tau_i, t) = \frac{\partial Q_i}{\partial t} + \nabla_{x_i} Q_i(x_i, t)\trn F_i(x_i, \tau_i, t),
$$
$\Gamma_i = \Gamma_i\trn > 0,$ $m_i\times m_i$ - matrix.

\begin{theorem}
\label{th_2_18}{\em
Suppose the
following assumptions hold for the system $\mathcal{S}$:
\begin{enumerate}
\item Functions $F_i(\cdot)$ are continuous in $x_i,$ $t,$ continuously
differentiable in $\tau_i$ and locally bounded in $t>0;$ functions $Q_i(x_i,t)$ are uniformly continuous in second argument
for all $x_i$ in bounded set, 
functions $\omega_i(x_i, \tau_i, t)$ are convex in $\tau_i;$
there exist constant vectors $\tau^*_i\in\mathbb{R}^{m_i}$ and scalar monotonically increasing
functions $\kappa_i(Q_i),\rho_i(Q_i)$ such that $\kappa_i(0)=\rho_i(0)=0, \lim_{Q_i\to+\infty}\kappa_i(Q_i)=+\infty$
\begin{equation}
\label{omega_leq_rho_2_18}
\omega_i(x_i, \tau^*_i, t)\leq -\rho_i(Q_i(x_i, t)), 
\end{equation}
and $Q_i(x_i,t)\geq\kappa_i(\|x_i-x_i^*(t)\|).$
\item functions $h_i(x, \tau, t)$ are continuous and satisfy the following inequalities
\begin{equation}
\label{grad_leq_rho_2_18}
|\nabla_{x_i}Q_i(x_i, t)\trn h_i(x, \tau, t)|\leq\sum_{j=1}^d \mu_{ij}\rho_j(Q_j(x_j, t)),
\end{equation}
where matrix $M-I$ is Hurwitz, $M=\{\mu_{ij}\},$ $\mu_{ij}>0,$ $I$ is identity matrix.
\end{enumerate}
Then system \eqref{obj_equat_2_18},\eqref{asg_2_18} is globally asymptotically stable in variables
$x_i-x_i^*(t),$ all trajectories are bounded on $t\in[0,+\infty)$ and satisfy \eqref{goal_funcs_2_18}.
}
\end{theorem}

\section{Main result}
\subsection{Problem statement. Adaptive controller structure}

Let the leader subsystem be described by the equation
\begin{equation}
\label{lead_eq}
\dot{\overline{x}} = A_L\overline{x}+B_L(\overline{u}+\psi_0(\overline{y})), \quad\overline{y}=C\trn \overline{x},
\end{equation}
where $\overline{x}\in\mathbb{R}^n$ -- state, $\overline{y}\in\mathbb{R}^l$ -- measurement, 
$\overline{u}(t)\in\mathbb{R}^1$ is control that specified in advance, 
$\psi_0\colon\mathbb{R}^l\to\mathbb{R}^1$ -- internal nonlinearity. Let $A_L, B_L, C$ and 
$\psi_0(\cdot)$ be known and not depending on the vector of unknown parameters $\xi\in\Xi,$ 
where $\Xi$ is known set.

Consider a network $S$ of $d$ interconnected subsystems
$S_i,\, i=1,\ldots, d, d\in\mathbb{N}.$ Let subsystem $S_i$ be described by the following
equation
\begin{equation}
\label{obj_eq}
\begin{aligned}
\dot{x}_i &= A_i x_i+B_i u_i+B_L\psi_0(y_i)+\sum_{j=1}^d\alpha_{ij}\varphi_{ij}(x_i-x_j),\\ 
y_i&=C\trn x_i,\qquad i=1,\ldots,d,\\
\end{aligned}
\end{equation}
where $x_i\in\mathbb{R}^n, u_i\in\mathbb{R}^1,
\alpha_{ij}\in\mathbb{R}^1, y_i\in\mathbb{R}^l.$ Functions
$\varphi_{ij}(\cdot),$ $i=1,\ldots,d,\,j=1,\ldots,d,$ describe
interconnections between subsystems. We assume $\varphi_{ii}=(0,0,\ldots,0)\trn, i=1,\ldots,d.$
Let matrices $A_i,B_i$ and 
functions $\varphi_{ij}(\cdot),$ $i=1,\ldots,d,j=1,\ldots,d,$ 
depend on the vector of unknown parameters $\xi\in\Xi.$

Network model \eqref{obj_eq} can describe, for example, interconnected electrical generators~\cite{Willems74}.

Let the control goal be specified as convergence of all
subsystems and the leader trajectories:
\begin{equation}
\label{goal}
\lim_{t\to+\infty}\left(x_i(t)-\overline{x}(t)\right) = 0,\qquad i=1,\ldots,d.
\end{equation}

The adaptive synchronization problem is to find a decentralized controller 
\\$u_i=\mathcal{U}_i(y_i, \overline{u}, t)$ 
ensuring the goal \eqref{goal} for all values of unknown plant parameters.

Denote $\sigma_i(t)=\mathrm{col}(y_i(t),\overline{u}(t)).$ Let the main loop 
of the adaptive system be specified as set of linear tunable local control laws:
\begin{equation}
\label{loc_contr}
u_i(t)=\tau_i(t)\trn\sigma_i(t),\quad i=1,\ldots,d,
\end{equation}
where $\tau_i(t)\in\mathbb{R}^{l+1}, i=1,\ldots,d$ are tunable parameters. By 
applying speed-gradient method~\cite{FMN99} it is easy to derive the following adaptation law:
\begin{equation}
\label{tau_contr}
\dot{\tau}_i=-g\trn (y_i-\overline{y})\Gamma_i\sigma_i(t), i=1,\ldots,d,
\end{equation}
where $\Gamma_i=\Gamma_i\trn >0$ -- $(l+1)\times(l+1)$ matrices, $g\in\mathbb{R}^l.$

\subsection{Synchronization conditions}

Introduce the following definition.

\begin{definition} 
\label{G_monot}
{\em Let $G\in\mathbb{R}^l.$ Function $f\colon\mathbb{R}^l\to\mathbb{R}^1$ is 
called G-monotonically decreasing if inequality 
$(x-y)\trn G \left(f(x)-f(y)\right)\leq 0$ holds for all $x,y\in\mathbb{R}^l$.}
\end{definition}

{\it Remark 1.} Apparently, for $l=1, G=1$ \quad $G$-monotonical decrease of
the function $f$ is equivalent to incremental passivity~\cite{PavlovMarconi}
of the static system with characteristics $(-f).$ Definition \ref{G_monot}
is easily extended to dynamical systems with the state vector 
$x\in\mathbb{R}^n,$ input $u\in\mathbb{R}^m$ and output $y\in\mathbb{R}^l.$
It corresponds to existence of a smooth function $V(x_1,x_2)$ satisfying
an integral inequality
$$
\dot{V}(x_1,x_2)\leq (u_2-u_1)\trn G(y_2-y_1).
$$
The corresponding property can be called incremental $G$-passivity by 
analogy with~\cite{Fradkov03}.

\par Consider real matrices $H=H\trn >0, g$ of size $n\times n, l\times 1$
correspondingly and a number $\rho>0$ such that:
\begin{equation}
\label{HA_L} 
HA_L+A_L\trn H <-\rho H,\quad HB_L=Cg.
\end{equation}
Denote $\lambda_*=\lambda_{max}(H)/\lambda_{min}(H)$ condition number of matrix $H,$
where $\lambda_{max}(H),\lambda_{min}(H)$ are maximum and minimum eigenvalues of matrix $H$.

For analysis of the system dynamics the following assumptions are made.

{A1)} The functions $\varphi_{ij}(\cdot),$$ i=1,\ldots,d,$ $j=1,\ldots,d$ are 
globally Lipschitz:
$$
\|\varphi_{ij}(x)-\varphi_{ij}(x')\|\leq L_{ij}\|x-x'\|,\quad L_{ij}>0. 
$$
The function $\psi_0(\cdot)$ is such that the unique existence of solutions of 
\eqref{lead_eq} holds.

{A2)}(Matching conditions,~\cite{SFEM}) For each $\xi\in\Xi$ there exist vectors $\nu_i=\nu_i(\xi)\in\mathbb{R}^l$ and
numbers $\theta_i=\theta_i(\xi)>0$ such that for $i=1,\ldots,d$
\begin{equation}
\label{compliance}
A_L=A_i+B_i\nu_i\trn C\trn, B_L=\theta_i B_i.
\end{equation}

Denote $\chi(s)=C\trn(sI_n-A_L)^{-1}B_L.$  For the case when matrix $A_L$ is Hurwitz introduce 
notation $\rho_*$ for stability degree of the function's $g\trn\chi(s)$ denominator, i.e. 
$\rho_*=\min_{k=1,\ldots,n}|\Re \lambda_{k}(A_L)|$ where $\lambda_{k}(A_L)$ are eigenvalues of $A_L$.

\begin{theorem}
\label{th_noident} 
{\em
Let $B_L\neq 0,$ matrix $A_L$ be Hurwitz and for some $g\in\mathbb{R}^l$ the following
frequency domain conditions hold:
\begin{equation}
\label{freq_compl}
\mathrm{Re}\, g\trn \chi(i\omega)>0,\quad \lim_{\omega\to\infty}\omega^2\,\mathrm{Re}\, g\trn\chi(i\omega)>0
\end{equation}
for all $\omega\in\mathbb{R}^1.$ Then there exist $H=H\trn>0, \rho>0$ such that relations \eqref{HA_L} hold.

Let for all $\xi\in\Xi$ Assumptions A1, A2 hold, function $\psi_0(\cdot)$ be $g$-monotonically decreasing, and
following inequalities hold
\begin{equation}
\label{connect}
\sum_{j=1}^d\left| \alpha_{ij}L_{ij}\right| <\gamma \quad i=1,\ldots,d,
\end{equation}
where $\gamma=\rho_*/(4d\lambda_*)$,
$\lambda_*$ is condition number of matrix $H$.

Then for all $\xi\in\Xi, i=1,\ldots,d$ adaptive controller \eqref{loc_contr},\eqref{tau_contr} ensures
achievement of the goal \eqref{goal} and boundedness of functions $\theta_i(t)$ on $[0,\infty)$ for all 
solutions of the closed-loop system \eqref{lead_eq}, \eqref{obj_eq}, \eqref{loc_contr}, \eqref{tau_contr}.
}
\end{theorem}

{\em Proof.}
Let's apply Lemma \ref{Yakub_lemma}. Note that in our case $m=1,$ 
i. e. $u$ is scalar. Let's choose $Cg$ instead of $C$ in \eqref{HA_Yak}. Then statement 
of the Lemma \ref{Yakub_lemma} and
conditions of Theorem \ref{th_noident} ensure existence of matrix
$H=H\trn>0$ such that 
$$HA_L+A_L\trn H<0,\quad HB_L=Cg.$$ Now we can conclude that there exists number $\rho>0$ such that 
the following is true:
\begin{equation}
\label{HA_rho}
HA_L+A_L\trn H<-\rho H,\quad HB_L=Cg.
\end{equation}

Denoting $z_i=x_i-\overline{x}$
introduce auxiliary error subsystems:
\begin{equation}
\label{temp_eq_noi}
\begin{aligned}
\dot{z}_i = &A_i x_i+B_i u_i+B_L\psi_0(y_i)+\sum_{j=1}^d\alpha_{ij}\varphi_{ij}(x_i-x_j)-
          \left(A_L\overline{x}+B_L(\overline{u}+\psi_0(\overline{y}))\right),\\
\tilde{y}_i=&C\trn z_i,\qquad i=1,\ldots,d,\\
\end{aligned}
\end{equation}
here we choose $u_i(t)$ same as in \eqref{loc_contr}.

Let us choose following goal functions 
$Q_i(z_i)=\frac12z_i\trn Hz_i,$  and apply Theorem \ref{th_2_18}. We need to 
evaluate the derivative trajectories of $Q_i(z_i)$ along trajectories of 
isolated (i.e. without interconnections)
auxiliary subsystems \eqref{temp_eq_noi}:
\begin{equation}
\label{omega_i}
\begin{aligned}
\omega_i(x_i,\overline{x},\tau_i) &= z_i\trn H[A_i x_i+B_i \tau_i\trn (t)\sigma_i(t)+
        B_L\psi_0(y_i)-A_L\overline{x}-B_L(\overline{u}+\psi_0(\overline{y}))].
\end{aligned}
\end{equation}
Denote $\tau_i^*=\mathrm{col}(\nu_i, \theta_i), i=1,\ldots,d.$ By taking $\tau_i=\tau_i^*,i=1,\ldots,d,$
we obtain
\begin{equation*}
\begin{aligned}
\omega_i(x_i,\overline{x},\tau_i^*)=&z_i\trn H[A_i x_i+B_i(\nu_i C\trn x_i+\theta_i\overline{u})+
B_L\psi_0(y_i)-A_L\overline{x}-B_L\overline{u}-B_L\psi_0(\overline{y})]=\\
&z_i\trn H[A_L x_i+B_L\overline{u}+B_L(\psi_0(y_i)-\psi_0(\overline{y}))-A_L x_i-B_L\overline{u}]=\\
&z_i\trn H[A_L z_i+B_L(\psi_0(y_i)-\psi_0(\overline{y}))].\\
\end{aligned}
\end{equation*}
Further, for $i=1,\ldots,d$
\begin{equation*}
z_i\trn HB_L(\psi_0(y_i)-\psi_0(\overline{y}))=z_i\trn Cg(\psi_0(y_i)-\psi_0(\overline{y}))=
(y_i-\overline{y})\trn g(\psi_0(y_i)-\psi_0(\overline{y}))\leq 0.
\end{equation*}
The last inequality holds because $\psi_0(\cdot)$ is $g$-monotonically decreasing. So
$$\omega_i(x_i,\overline{x},\tau_i^*)\leq \frac12 z_i\trn (H A_L+A_L\trn H) z_i.$$
Taking into account \eqref{HA_rho} we conclude 
$$
\omega_i(x_i,\overline{x},\tau_i^*)\leq -\rho\, Q_i(z_i).
$$
By taking $\rho_i(Q)=\rho\cdot Q$ we ensure that \eqref{omega_leq_rho_2_18} holds for $i=1,\ldots,d$. 
Other conditions from the first part of Theorem \ref{th_2_18} hold, since the
right hand side of the system \eqref{temp_eq_noi} and function $Q_i(z_i)$
are continuous in $z_i$ functions not depending in $t$ for any $i=1,\ldots,d$. Convexity condition is valid since the right hand side of 
\eqref{omega_i} is linear in $\tau_i$.

The interconnection condition \eqref{grad_leq_rho_2_18} in our case reads:
\begin{equation}
\label{grad_Q}
|\nabla_{z_i}Q(z_i)\trn \sum_{j=1}^d\alpha_{ij}\varphi_{ij}(z_i-z_j)|\leq\sum_{j=1}^d \mu_{ij}\rho\cdot Q(z_j),
\end{equation}
where $i=1,\ldots,d,$ and
 matrix $M-I$ should be Hurwitz ($M=\{\mu_{ij}\},$
$\mu_{ij}> 0$). 

For the case $d=1$ we can take $\mu_{11}=0.5$ and last inequality will be satisfied. Let's consider case $d>1.$

For $i=1,\ldots,d$ rewrite \eqref{grad_Q} as follows:
\begin{equation}
\label{ziH}
|z_i\trn H\sum_{j=1}^d\alpha_{ij}\varphi_{ij}(z_i-z_j)| \leq \frac{\rho}2\sum_{j=1}^d \mu_{ij}\, z_j\trn H z_j.
\end{equation}
Evaluate the left-hand side of \eqref{ziH}:
\begin{equation*}
\begin{split}
&\left|z_i\trn H\sum_{j=1}^d\alpha_{ij}\varphi_{ij}(z_i-z_j)\right|\leq
\sum_{j=1}^d\left|z_i\trn H\alpha_{ij}\varphi_{ij}(z_i-z_j)\right|\leq\\
&\sum_{j=1}^d|\alpha_{ij}L_{ij}|\cdot\|z_i\|\cdot\|H\|\cdot\|z_i-z_j\|\leq
\sum_{j=1}^d|\alpha_{ij}L_{ij}|\cdot\lambda_{max}(H)\cdot(\|z_i\|^2+\|z_i\|\cdot\|z_j\|),
\end{split}
\end{equation*}
for $i=1,\ldots,d.$ Then for $i=1,\ldots,d$ evaluate lower bound of the right-hand side of \eqref{ziH}:
$$
\frac{\rho}2\,\sum_{j=1}^d\mu_{ij} z_j\trn Hz_j \geq \frac{\rho}2\,\sum_{j=1}^d\mu_{ij} \lambda_{min}(H)\|z_j\|^2.
$$
It is seen that for $i=1,\ldots,d$ to ensure \eqref{grad_leq_rho_2_18} it is sufficient to impose an
inequality
$$
\sum_{j=1}^d|\alpha_{ij}L_{ij}|\cdot\lambda_{max}(H)\cdot(\|z_i\|^2+\|z_i\|\cdot\|z_j\|)\leq
\frac{\rho}2\,\sum_{j=1}^d\mu_{ij} \lambda_{min}(H)\|z_j\|^2,
$$
or
\begin{equation}
\label{temp_eq}
\sum_{j=1}^d|\alpha_{ij}L_{ij}|\cdot(\|z_i\|^2+\|z_i\|\cdot\|z_j\|)\leq
\frac{\rho}2\,\frac{\lambda_{min}(H)}{\lambda_{max}(H)}\sum_{j=1}^d\mu_{ij} \|z_j\|^2,
\qquad i=1,\ldots,d.
\end{equation}

Denote $\zeta=\rho/(\alpha_{max}\cdot 2\lambda_*),$ where
$$\alpha_{max}=\max_{i:1\leq i\leq d}\,\sum_{j=1}^d\left| \alpha_{ij}L_{ij}\right|.$$
Noting that $\rho$ in \eqref{HA_L} can be chosen arbitrarily close to $\rho_*$ and taking into account \eqref{connect} we can conclude that 
$$
\zeta>2d.
$$

The left-hand side of \eqref{temp_eq}:
\begin{equation*}
\begin{aligned}
\sum_{j=1}^d |\alpha_{ij}L_{ij}|(\|z_i\|^2+\|z_i\|\cdot &\|z_j\|)\leq 
\left(\sum_{j=1}^d|\alpha_{ij}L_{ij}|\right)\cdot\left(\sum_{j=1}^d(\|z_i\|^2+\|z_i\|\cdot\|z_j\|)\right)\leq\\
&\frac 12\,\alpha_{max}\left(3d\,\|z_i\|^2+\sum_{j=1}^d\|z_j\|^2\right),\quad i=1,\ldots,d.
\end{aligned}
\end{equation*}

Thus, if following inequality holds then \eqref{grad_leq_rho_2_18} is ensured:
\begin{equation}
\label{temp_eq2}
3d\,\|z_i\|^2+\sum_{j=1}^d\|z_j\|^2 \leq 2\zeta\cdot\sum_{j=1}^d \mu_{ij} \|z_j\|^2,\quad i=1,\ldots,d.
\end{equation}

Introduce matrix $M=\{\mu_{ij}\}$ as follows
$$M=\begin{pmatrix}
\mu_{11}&\mu_{12}&\ldots&\mu_{1d}\\
\mu_{21}&\mu_{22}&\ldots&\mu_{2d}\\
\vdots&\vdots&\ddots&\vdots\\
\mu_{d1}&\mu_{d2}&\ldots&\mu_{dd}
\end{pmatrix},
\qquad\mu_{ij}=\left\{
\begin{aligned}
&\frac1{2\zeta}(3d+1),&\quad i=j;\\
&\frac1{2\zeta},&\quad i\neq j.\\
\end{aligned}
\right.
$$
Such choice of $M$ ensures \eqref{temp_eq2}.

Note that $M$ is symmetric. If matrix $I-M$ is positive definite then $M-I$ is Hurwitz. 
Diagonal elements of $I-M$ are positive since $d>1$ and $\zeta>2d.$ By taking into account that
$$
1-\frac{1}{2\zeta}(3d+1)-(d-1)\frac{1}{2\zeta}>0
$$
and applying Gershgorin circle Theorem we conclude that $M-I$ is positive definite.

Thus, statement of the Theorem \ref{th_noident}
follows from Theorem \ref{th_2_18}.
\square

{\it Remark 2.} The value of $\gamma$ can be evaluated by solving 
LMI \eqref{HA_rho} by means of one of existing software package.

{\it Remark 3.} By interconnections graph of network $S$ we can consider directed 
graph which is a pair of two sets: a set of nodes and a set of arcs. Cardinality
of a set of nodes is $d;$ $i$-th node is associated with subsystem $S_i$ for any
$i=1,\ldots,d.$ We say that arc from $i$-th node to $j$-th node belongs to the set of 
arcs if $\varphi_{ij}(\cdot)$ is not zero function. By weighted in-degree of $i$-th
node we define following number: $\sum_{j=1}^d\left| \alpha_{ij}L_{ij}\right|.$
If each nonzero addend from last sum is equal to 1 then introduced definition of weighted 
in-degree of the node coincides with the definition of in-degree of digraph's node. 
Thus the inequality \eqref{connect} can be interpreted as follows: weighted in-degree
of each node of interconnections graph must be less than $\gamma.$

\section{Example. Network of Chua circuits}
\subsection{System description and theoretical study}
Chua circuit is a well known example of simple nonlinear system
possessing complex chaotic behavior \cite{WuChua95}. Its
trajectories are unstable and it is represented in
the Lurie form. Let us apply our results to synchronization 
with leader subsystem in
the network of five interconnected nonidentical Chua systems.

Let $m_0=-8/7, m_1=-5/7, p=15.6, q=30, b=1$ and $g=1.$

Let the leader subsystem be described by the equation
$$\dot{\overline{x}} = A_L\overline{x}+B_L(\overline{u}+\psi_0(\overline{y})), \quad\overline{y}=C\trn \overline{x},$$
where $\overline{x}\in\mathbb{R}^3$ is state vector of the system, $\overline{y}\in\mathbb{R}^1$ is output 
available for measurement, $\overline{u}$ is scalar control variable,
$\psi_0(\overline{y})=pv(\overline{y})/b,$ where $v(x)=-0.5(m_0-m_1)(|x+1|-|x-1|-2x).$ 
Further, let
$$
A_L=
\begin{pmatrix}
-1&0&0\\
1&-1&1\\
0&-q&0
\end{pmatrix},
$$
$B_L=\mathrm{col}(b,0,0), C=\mathrm{col}(1,0,0).$ 

\begin{figure}
\centering
\includegraphics[width=3in,height=1.8in]{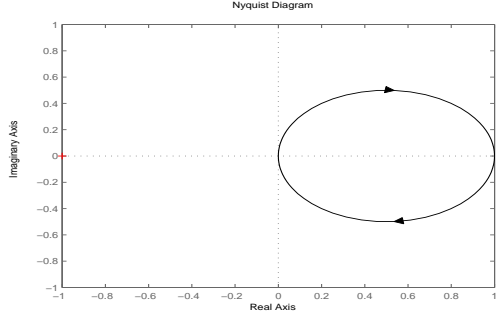}
\caption{Nyquist plot of $\chi(i\omega), \omega\in\mathbb{R}^1$.}
\label{nyquist}
\vspace*{2pt}
\end{figure}

Transfer function $\chi(s)=C\trn(sI-A_L)^{-1}B_L=(s^2+s+30)/(s^3+2s^2+31s+30).$ 
It is seen from the Nyquist plot of $\chi(i\omega),\forall \omega\in\mathbb{R}^1,$
presented on Fig. \ref{nyquist}, that first frequency domain inequality of \eqref{freq_compl}
holds. The second frequency domain inequality of \eqref{freq_compl} also holds since relative
degree of $\chi(s)$ is equal to one and highest coefficient of its numerator is positive.

Obviously $\psi_0(\cdot)$ is $g$-monotonically decreasing. 

Let subsystem $S_i$ for $i=1,\ldots,5$ be described by \eqref{obj_eq} with $u_i,\alpha_{ij}\in\mathbb{R}^1.$
By choosing $(\nu_1,\nu_2,\nu_3,\nu_4,\nu_5)=(3,1,4,1,5), \theta_i=1/i, i=1,\ldots,5$ 
and using \eqref{compliance} we obtain matrices $A_i,B_i$ for $i=1,\ldots,5,$ which are not equal, i.e.
nodes are nonidentical.
Denote 
$\varphi_{ij}=\varphi_{ij}(x_i-x_j), i=1,\ldots,5, j=1,\ldots,5.$
Let $\varphi_{14}, \varphi_{25}, \varphi_{32}, \varphi_{42}, \varphi_{45}, \varphi_{52}, \varphi_{53},$ 
be equal to $(0,0,0)\trn.$ Further, let 
\begin{equation*}
\begin{aligned}
&\varphi_{12}=(\sin(x_{11}-x_{21}),0,0)\trn,\qquad &\varphi_{13}&=(0,x_{12}-x_{32},0)\trn,\\
&\varphi_{15}=(0,0,\sin(x_{13}-x_{53}))\trn,\qquad &\varphi_{21}&=(x_{21}-x_{11},0,x_{23}-x_{13})\trn,\\
&\varphi_{23}=(0,\sin(x_{22}-x_{32}),0)\trn,\qquad &\varphi_{24}&=(0,x_{22}-x_{42},0)\trn,\\
&\varphi_{31}=(\sin(x_{31}-x_{11}),0,0)\trn,\qquad &\varphi_{34}&=(\sin(x_{31}-x_{41}),0,0)\trn,\\
&\varphi_{35}=(x_{31}-x_{51},x_{32}-x_{52},x_{33}-x_{53})\trn,\qquad &\varphi_{41}&=(0,\sin(x_{42}-x_{12}),0)\trn,\\
&\varphi_{43}=(\sin(x_{41}-x_{31}),0,0)\trn,\qquad &\varphi_{51}&=(x_{51}-x_{11},0,x_{53}-x_{13})\trn,\\
&\varphi_{54}=(0,x_{52}-x_{42},0)\trn.&&\\
\end{aligned}
\end{equation*}
Lipschitz constants of all $\varphi_{ij}$ are equal to $1.$

It follows from Theorem \ref{th_noident} that decentralized adaptive control \eqref{loc_contr} 
provides synchronization goal \eqref{goal}
if for all $i=1,\ldots,5$ inequality $\sum_{j=1}^5\left| \alpha_{ij}\right|<\gamma$ holds,
i.e. if interconnections are sufficiently weak.

\subsection{Simulation results}

Consider following control of leader subsystem 
$\overline{u}=\frac 1b\left[(-(1+m_0)p+1)\overline{x}_1+p\overline{x}_2\right].$
Such $\overline{u}$ ensures chaotic behavior of leader subsystem.
Let us put $\Gamma_i = I, i=1,\ldots,d,$ where $I$ -- identity matrix, and
\begin{equation*}
\begin{aligned}
\overline{x}_1(0) &= 0.5,\quad\overline{x}_2(0) = 0,\quad \overline{x}_3(0) = 0,\\
x_1(0)&=(7,14,0.4)\trn, \quad x_2(0)=(0,4,4)\trn \\
x_3(0)&=(1,-1,4.5)\trn, \quad x_4(0)=(3,-4,0.2)\trn \\
x_5(0)&=(2,8,15).
\end{aligned}
\end{equation*}

\begin{figure*}[!t]
\centering
\includegraphics[width=5.6in,height=2.6in]{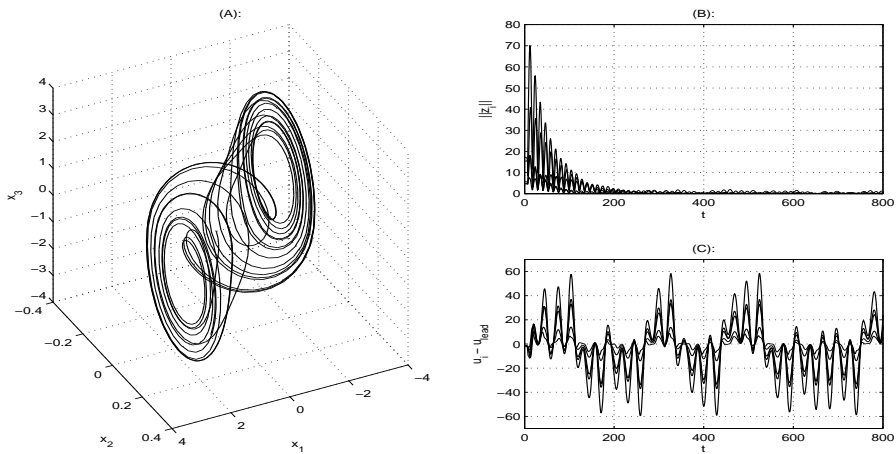}
\caption{(A): Phase portrait of leader subsystem,\quad (B): $\|z_i\|$,\quad 
   (C): $\tilde{u}_i=u_i-\overline{u}, i=1,\ldots,5.$}
\label{fig_sim}
\vspace*{2pt}
\end{figure*}

Denote by $\alpha$ $5\times 5$ matrix with element $\alpha_{ij}$ lying in
the $i$-th row and the $j$-th column, $i,j=1,\ldots,5,$ and
$$
\widehat{\alpha}=\begin{pmatrix}
     0& 0.0051&   0.1395&         0&    0.1676\\
0.0662&     0&   0.0921&    0.0065&         0\\
0.2013&     0&        0&    0.2271&    0.1430\\
0.0907&        0&   0.0675&         0&         0\\
0.0663&        0&        0&    0.2773&         0\\
\end{pmatrix}.
$$
Let us choose adaptive control $u_i, i=1,\ldots,5$ as in \eqref{loc_contr} 
and apply Theorem \ref{th_noident}.
If we take $\alpha=\widehat\alpha,$ then simulation shows that 
$\|z_i\|\to 0, i=1,\ldots,5,$ i.e. synchronization is achieved:
all state vectors of nonidentical nodes converge to the state vector of the leader subsystem,
see Fig. \ref{fig_sim}-(B).
Phase portrait of the leader subsystem, $\|z_i\|, \tilde{u}=u_i-\overline{u}, i=1,\ldots,5$ 
found by 40 sec. simulation are shown on Fig.
\ref{fig_sim}.

\section{Conclusions}

In contrast to a large number of previous results, we obtained
synchronization conditions for networks consisting of nonidentical
nonlinear systems with
incomplete measurement, incomplete control, incomplete information 
about system parameters and coupling. The
design of the control algorithm providing synchronization
property is based on speed-gradient method
\cite{FMN99}, while derivation of synchronizability
conditions is based on Yakubovich-Kalman lemma and result presented 
in \cite{Fradkov_1990}.

\end{document}